\documentclass[a4paper, 12pt]{article}

\usepackage{amsmath, bm}
\usepackage{amsfonts}
\usepackage{amsthm}
\usepackage{enumitem}
\usepackage{ascmac}               
\usepackage{multicol}             
\usepackage[pdftex]{graphicx}     
\usepackage{float}                
\usepackage{wrapfig}              
\usepackage[hang,small,bf]{caption}
\usepackage[subrefformat=parens]{subcaption}
\usepackage{here}

\usepackage{indentfirst} 

\theoremstyle{plain}

\theoremstyle{definition}
\newtheorem{dfn}{Def}
\newtheorem{thm}[dfn]{Thm}

\newtheorem{exam}[dfn]{Example}
\newtheorem*{rem}{Remark}

\title{Cluster Analysis for Globally Coupled Map
using Optimal Transport Distance \\
and Complexity of Attractor-ruin}
\author{
	Koji Wada
	\thanks{Department of Mathematics, Graduate School of Science, Hokkaido University,
	Kita 10, Nishi 8, Kita-Ku, Sapporo, Hokkaido, 060-0810, Japan.
	wada.koji.n4@elms.hokudai.ac.jp}
	\quad and \quad
	Takao Namiki
	\thanks{Department of Mathematics, Faculty of Science, Hokkaido University,
	Kita 10, Nishi 8, Kita-Ku, Sapporo, Hokkaido, 060-0810, Japan.
	nami@math.sci.hokudai.ac.jp}
	}
\date{}

\begin{document}
\maketitle

\begin{description}[leftmargin=13ex]
	\item[\bf Keywords:] Dynamical system, Glocally coupled map, Chaotic itinerancy, \\
	Optimal transport distance, Attractor ruin,
	\item[{\bf MSC2020:}] 37D45, 39A33, 62H30
\end{description}

\begin{abstract}
  In this paper,
  we show the results of the strength of attractor-ruins for a globally coupled map.
  The globally coupled map (GCM) is a discrete dynamical system,
  and here we consider a model in which the logistic map is globally coupled.
  An attractor-ruin is a set in which the attractor is destabilized by a change in parameters,
  which is characterized by a Milnor attractor.
  Intermittent phenomena called chaotic itinerancy,
  in which orbits transition between attractor-ruin,
  have been observed in various complex systems,
  and their onset mechanisms and statistical properties have attracted attention.
  In this study,
  the instability of orbits of GCM is analyzed from the perspective of clustering using the optimal transport distance,
  and the strength of attractor-ruins is numerically evaluated by applying this method.
  As a result,
  it was found that the strength of various attractor-ruins is high in the parameter region called the partially ordered phase,
  where chaotic itinerancy occurs.
\end{abstract}

\section{Introduction}\label{Intro}
Nonlinear phenomena called ``chaos'' are observable in a one-dimensional discrete dynamical system such as a logistic map,
and several mathematical definitions have been proposed by Li and Yorke \cite{li1975period}, R.Devaney \cite{devaney1986introduction}, and others.
Since the 1980s, in low-dimensional dynamical systems, the characteristics of nonlinear phenomena called “intermittency” have been investigated \cite{pomeau1980intermittent, grebogi1987critical, ott1994blowout}.
Then, an intermittent phenomenon,
which is called Chaotic itinerancy (CI),
is independently discovered
in a model of optimal turbulence \cite{ikeda1989maxwell},
in a globally coupled system \cite{kaneko1990clustering},
and in nonequilibrium neural networks \cite{tsuda1991chaotic}.
CI is the phenomenon transition between a chaotic set, called  attractor-ruin, in which the attractor is destabilized by a change in parameters, and an attractor-ruin is mathematically characterized by a Milnor attractor \cite{milnor1985concept, kaneko1997dominance}.
CI is considered a universal phenomenon in a class of high-dimensional dynamical systems [9], and its
statistical properties and mechanisms are being actively investigated,
but without sufficient results to understand the phenomenon.

In this paper,
we focus on a mathematical model called the globally coupled map (GCM) proposed by K.Kaneko \cite{kaneko1990clustering}.
GCM is known to have invariant sets that can be interpreted in perspective of clustering (see Sec.~I\hspace{-1.2pt}I) and are known to be attractors for many parameters.
Then, for some parameter regions,
clusters are unstable and GCM can be interpreted as generating CI.
However, no attempt has been made to quantitatively evaluate such changes in clustering patterns to capture CI precisely.
Therefore, in this study,
we evaluated the changes in the clustering pattern using the optimal transport distance and information on the number of elements in each cluster.

The outline of this paper is as follows:
First, Sec.~I\hspace{-1.2pt}I briefly introduces features of orbits for GCM in each parameter region.
GCM is an $n$-dimensional discrete dynamical system
and here we consider a coupled system of logistic maps $f(x) = \alpha x(1- x)\ (x, \alpha \in \mathbb{R})$.
In Sec.~I\hspace{-1.2pt}I\hspace{-1.2pt}I,
a definition of an optimal transport distance is introduced
and a definition of cluster distributions is given in preparation for the analysis.
In Sec.~I\hspace{-1.2pt}V,
the main results of numerical studies are presented.
We observed intermittent phenomena in the partially ordered phase (Fig.~\ref{fig:3-c})
and obtained a higher time averages of the time series of optimal transport distance (Fig.~\ref{fig:5}).
In addition,
we evaluated the strength of attractor-ruins using the optimal transport distance and the effective dimension.
As a result, the strength of attractor-ruins is higher in the partially ordered phase.
Finally, in Sec.~V, a summary and discussions are presented.

\section{Globally coupled map and its features}
The globally coupled map (GCM) is proposed by K.Kaneko\cite{kaneko1990clustering} as a mathematical model in which CI is observed.
GCM is a $N$-dimensional discrete dynamical system,
which globally coupled nonlinear maps.
So, we consider here a coupled system of logistic maps,
\begin{gather*}
	x_{n+1}(i) = (1 - \varepsilon)f(x_n(i)) + \frac{\varepsilon}{N}\sum_{j=1}^Nf(x_n(j)),
	\quad
	i = 1,2,\ldots,N\\
	f(x) = \alpha x(1 - x),
	\quad
  \bm{x}_0 \in [0, 1]^N,
  \quad
	(\alpha, \varepsilon) \in [0, 4]\times [0, 1]
	\end{gather*}
where $n$ is discrete time
and $i$ is the index of elements.
Note that the orbits do not diverge
because the domain of GCM is considered on $[0,1]^N$.
GCM has two parameters which are called a bifurcation parameter $\alpha$ and a coupling coefficient $\varepsilon$.
Depending on these parameters $\alpha$ and $\varepsilon$,
GCM forms clusters of various types in which each element oscillates synchronously $(x(i) = x(j),\ i\neq j)$.
Convergence of a cluster concerning GCM can be interpreted as an attractor to some invariant set $H_\sigma$ in the phase space \cite{komuro1999kokyuroku}.

\begin{figure}[t]
	\centering
	\begin{minipage}[t]{0.45\columnwidth}
		\centering
		\includegraphics[width=0.9\columnwidth]{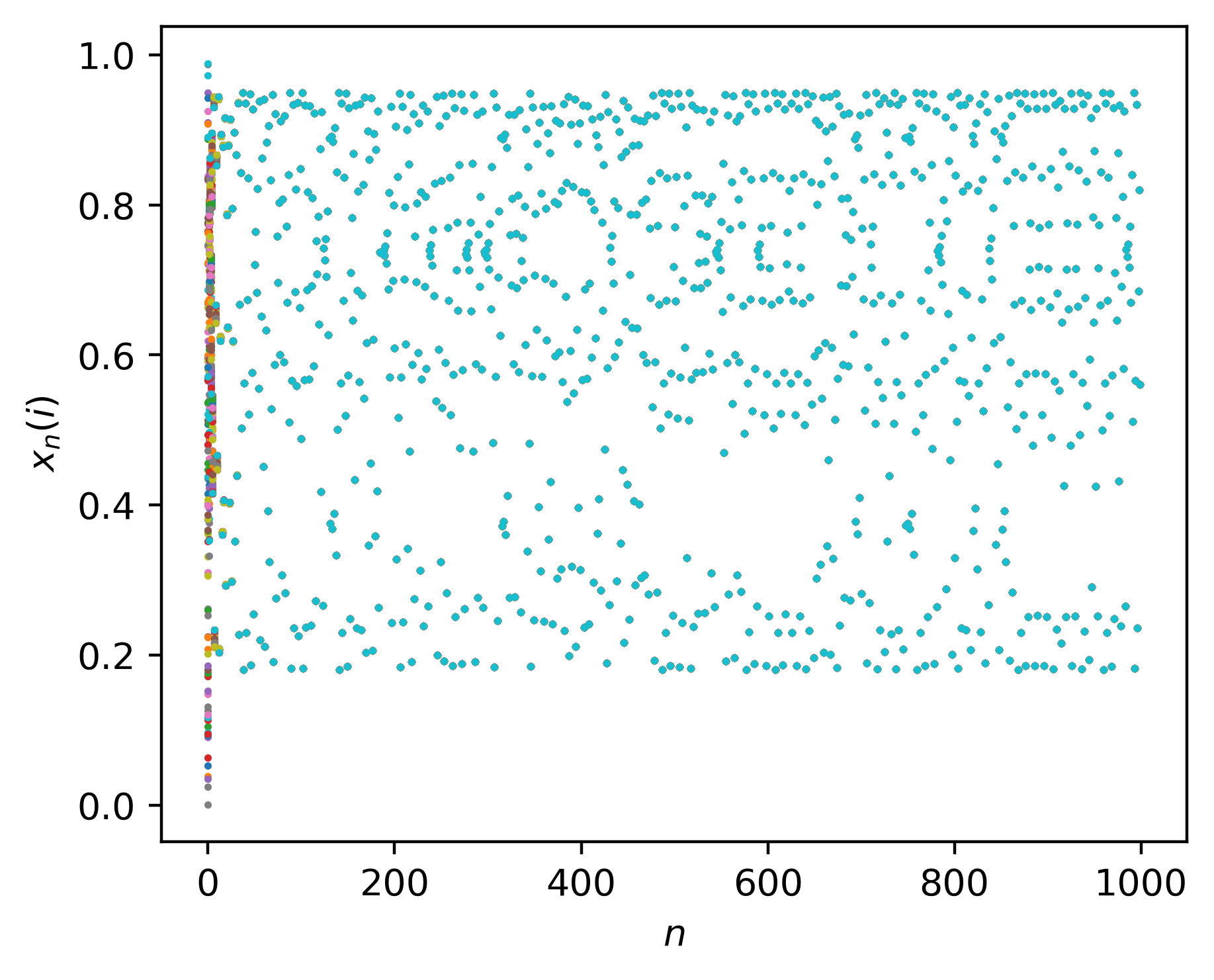}
		\subcaption{Coherent phase $(\varepsilon = 0.5)$}
	\end{minipage}
	\hspace{0.05\columnwidth}
	\begin{minipage}[t]{0.45\columnwidth}
		\centering
		\includegraphics[width=0.9\columnwidth]{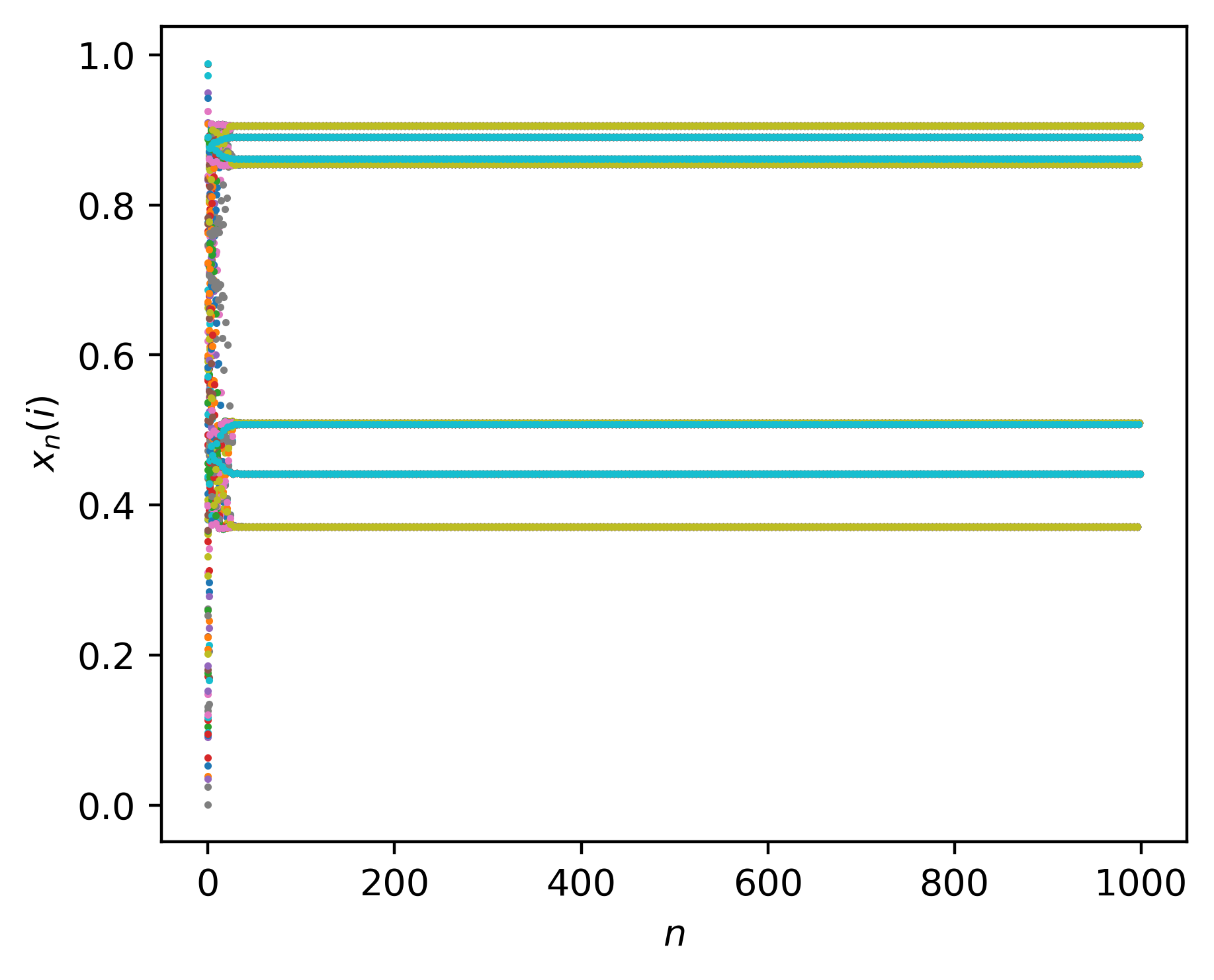}
		\subcaption{Ordered phase $(\varepsilon = 0.2)$}
	\end{minipage}
	\vspace{10pt}\\
	\begin{minipage}[b]{0.45\columnwidth}
		\centering
		\includegraphics[width=0.9\columnwidth]{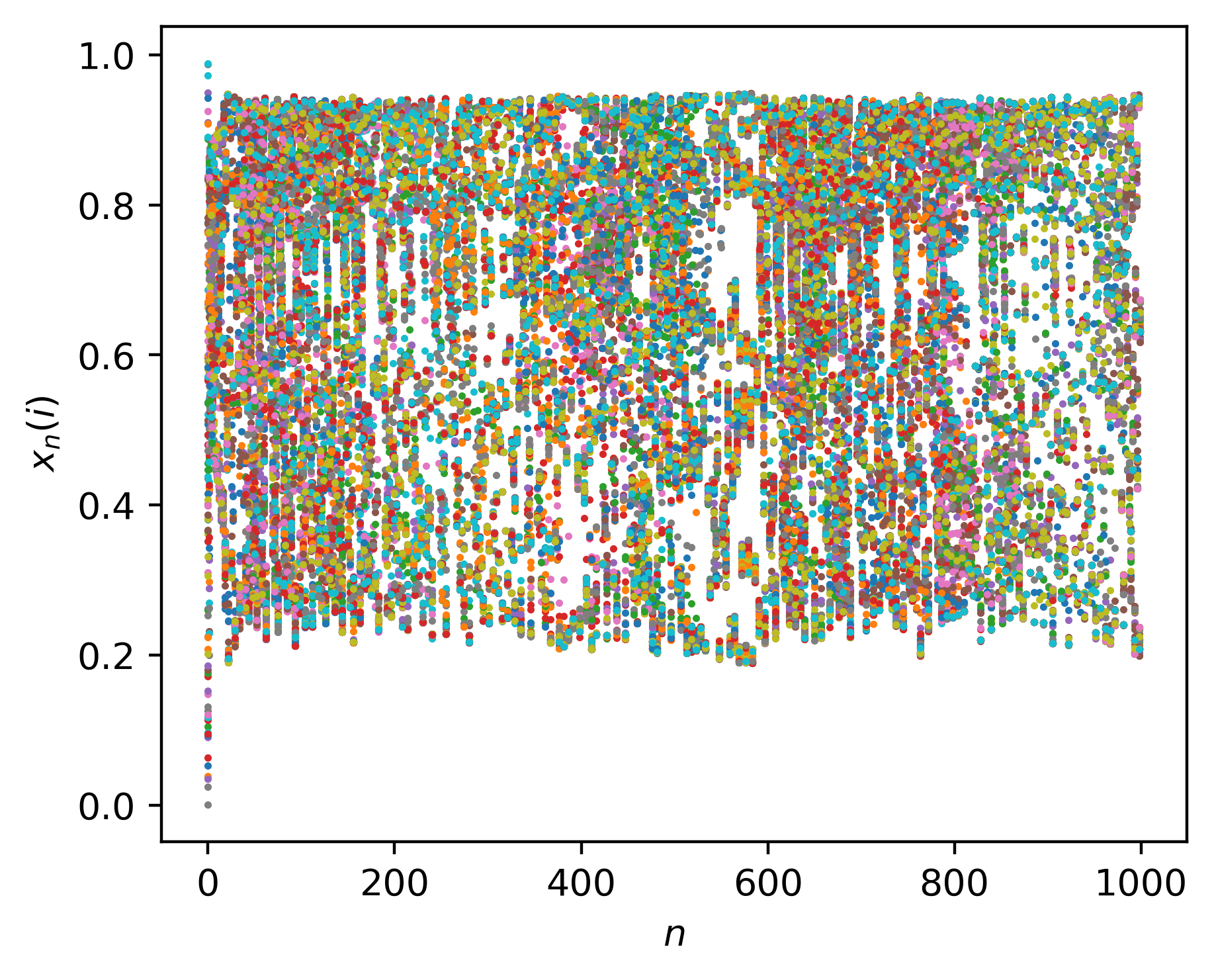}
		\subcaption{Partially ordered phase $(\varepsilon = 0.3)$}
	\end{minipage}
	\hspace{0.05\columnwidth}
	\begin{minipage}[b]{0.45\columnwidth}
		\centering
		\includegraphics[width=0.9\columnwidth]{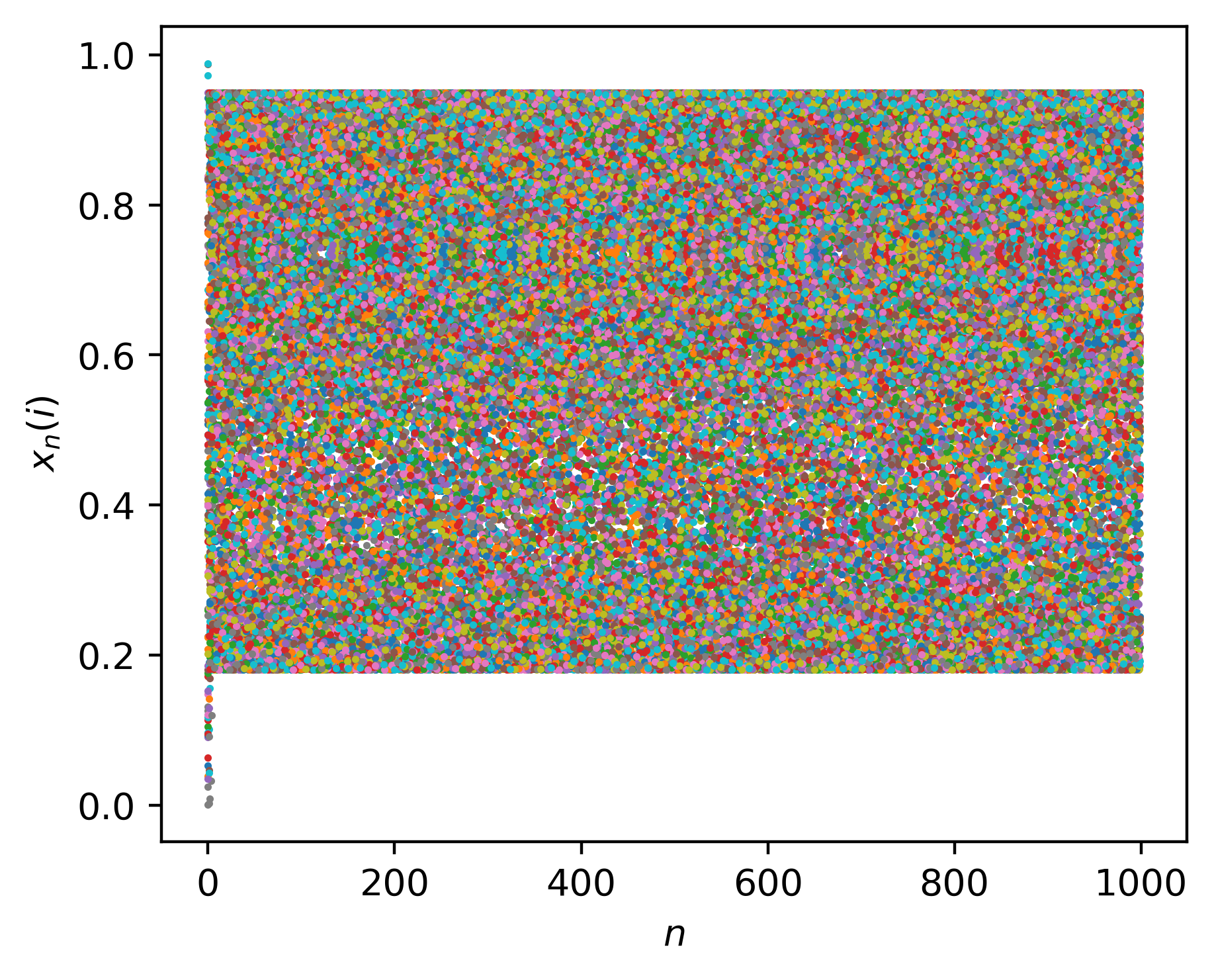}
		\subcaption{Turbulent phase $(\varepsilon = 0.0)$}
	\end{minipage}
	\caption{Examples of GCM orbits for each phase ($N = 100, \alpha = 3.8$).
	The initial values of these orbits are identical.}
	\label{fig:1}
\end{figure}

\begin{dfn}\label{dfn:effdim}
	Let $S_N$ be a symmetry group of degree $N$.
	For any $\sigma \in S_N$,
	\begin{gather*}
		H_\sigma := \{\bm{x} = (x(1), x(2), \ldots, x(N)) \in [0,1]^N \mid x(i) = x(\sigma(i)),\ i =1,2,\ldots,N\}.
	\end{gather*}
	In addition,
	for any $x \in [0,1],\ \delta > 0$,
	\begin{gather*}
		\textrm{ED}(x, \delta) := \inf_{\sigma \in S_N}\{\dim H_\sigma \mid H_\sigma \cap B_\delta(x) \neq \emptyset\}
	\end{gather*}
	is called the effective dimension.
\end{dfn}
It is known that in GCM,
the invariant set that serves as an attractor changes depending on the parameters $(\alpha, \varepsilon)$,
and in some cases,
clusters do not converge but continue to change unstably.
Based on characteristics of the evolution of clustering patterns,
the parameters can be classified into the following four categories
(Fig.~\ref{fig:1} shows examples of orbits for each phase in GCM) :
\begin{enumerate}
	\item Coherent phase: A parameter region where all elements $x(i)$ synchronize and converge into a single cluster,
	primarily found in the region where $\alpha$ is small and $\varepsilon$ is large.
	\item Ordered phase: A parameter region where the system converges into a small number of clusters.
	After convergence, GCM behaves as a system with a few degrees of freedom.
	\item Partially ordered phase: A parameter region where the effective dimension does not converge,
	and the orbit transitions through the neighborhoods of various invariant sets.
	In particular, this partially ordered phase is further classified into two types,
	Partially Ordered Phase I and Partially Ordered Phase I\hspace{-1.2pt}I.
	\item Turbulent phase: A parameter region where all elements $x(i)$ behave independently without synchronization,
	typically found in the region where $\alpha$ is large and $\varepsilon$ is small.
\end{enumerate}

\begin{figure}[t]
	\centering
	\includegraphics[width=0.7\columnwidth]{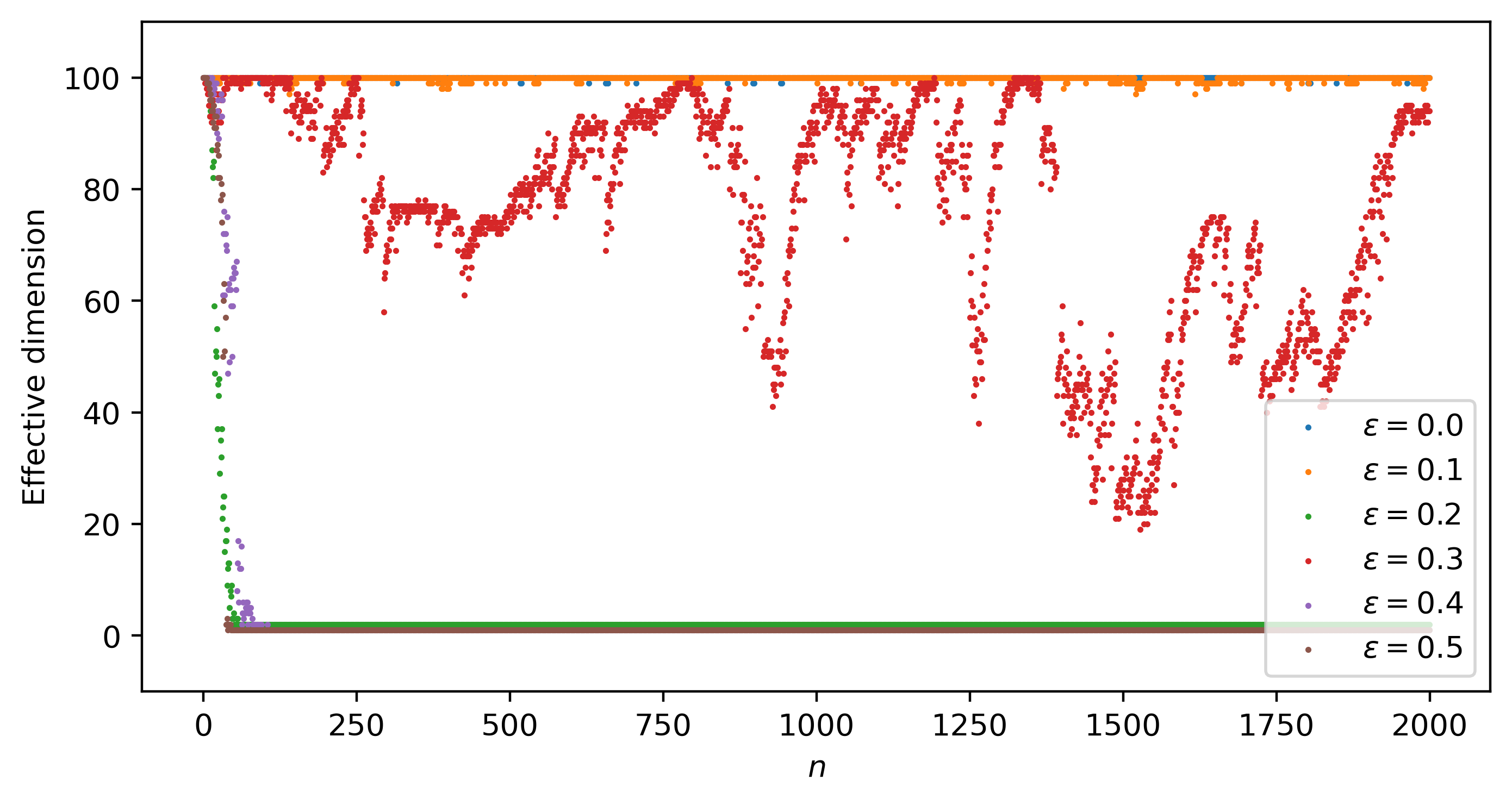}
    \caption{Examples of changes in the effective dimension for each parameter $\varepsilon$.
	Here, $N = 100$, $\alpha = 3.8$, $\delta = 10^{-6}$.
	The initial conditions are identical for all parameter settings.}
	\label{fig:2}
\end{figure}

Fig.~\ref{fig:2} shows the time series of the effective dimension when $\alpha = 0.3$ is fixed and $\varepsilon$ is varied.
The parameters $(\alpha, \varepsilon)$ corresponding to the coherent phase,
ordered phase, and turbulent phase quickly reach a steady state.
However, in the partially ordered phase, e.g., $\varepsilon = 0.3$,
the effective dimension does not converge and fluctuates unstably.

\begin{figure}[t]
	\centering
	\includegraphics[width=0.7\columnwidth]{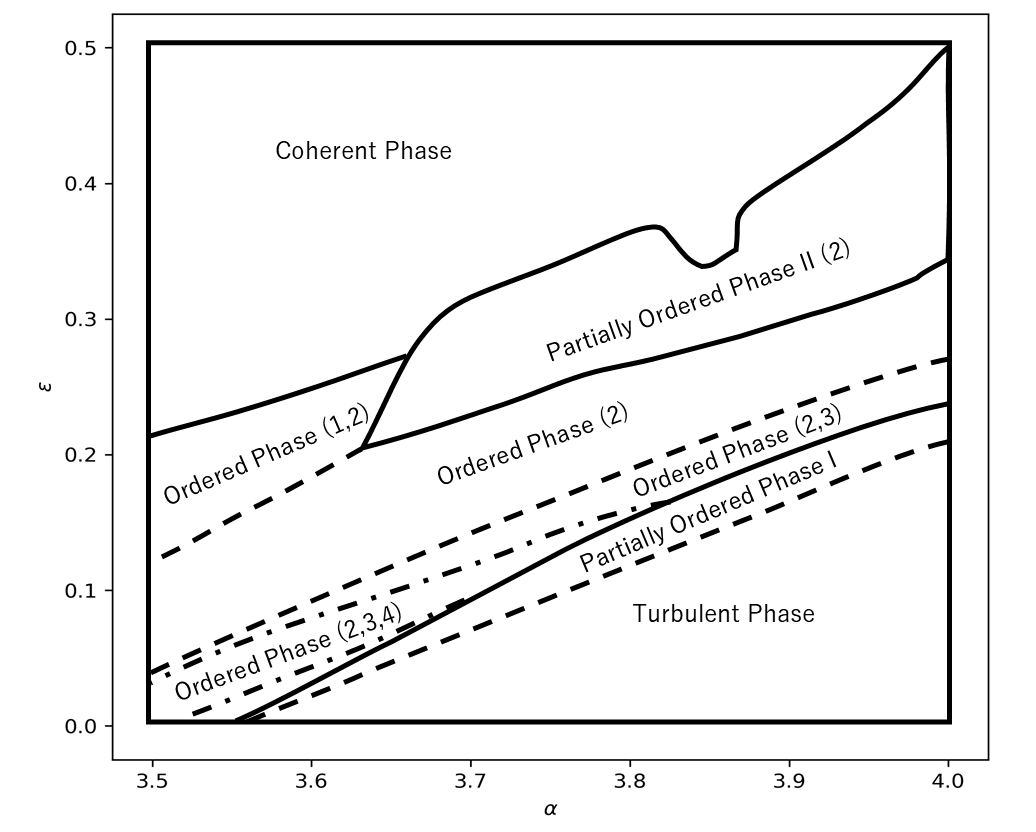}
    \caption{The phase diagram of GCM.}
	\label{fig:phase_diagram}
\end{figure}

Fig.~\ref{fig:phase_diagram} shows the phase diagram of GCM.
The value in parentheses within the figure represents the effective dimension (number of clusters) to which the system converges when the initial condition is appropriately set.

\section{Optimal transport distance and cluster distribution}
The optimal transport distance is the solution to the optimal transport problem,
which seeks to move a sandpile from one location to another at minimal cost.
This problem was first proposed by G.Monge in 1781
and was later given a well-defined formulation by L.Kantrovich \cite{villani2008optimal}.
The optimal transport problem can be succinctly described as a variational problem in the space of probability measures.
In particular, the optimal transport distance, known as the Wasserstein distance, is known to satisfy the axioms of a metric.

\begin{dfn}[\cite{villani2021topics}]
	Let $\mu$ and $\nu$ be probability measures defined on the measurable space $X$ and $Y$,
	respectively, and let $c : X \times Y \to [0, \infty]$ be a measurable function.
	Furthermore, consider a probability measure $\pi$ on $X \times Y$ that satisfies the following condition for any measurable sets $A \subset X$ and $B \subset Y$ :
	\begin{gather*}
		\pi(A \times Y) = \mu(A),
		\quad
		\pi(X \times A) = \nu(A).
	\end{gather*}
	Let $\Pi(\mu, \nu)$ denote the set of all such probability measures $\pi$.
	Here, $c$ is called the cost function,
	and $\pi$ is referred to as a transference plan.
	The optimal transportation cost between $\mu$ and $\nu$ is then defined as :
	\begin{gather*}
		T_c(\mu, \nu) := \inf_{\pi \in \Pi(\mu, \nu)}\int_{X \times Y}c(x, y)d\pi(x, y),
	\end{gather*}
	where the optimal $pi$ achieving this infimum is called the optimal transference plans.
\end{dfn}

By imposing appropriate assumptions on the cost function $c$,
it is possible to define a metric on the space of probability measures.

\begin{dfn}
	Let $X$ be a Polish space endowed with a distance $d$ and $p \in [1, \infty)$.
	We consider the cost function $c(x, y) = d(x, y)^p$ and $P_p(X)$,
	which is the set of probability measures $\mu$ such that for some (and thus any) $x_0 \in X$,
	\begin{gather*}
		\int d(x_0, x)^p d\mu(x) < \infty.
	\end{gather*}
	$W_p := T_d^{1/p}$ is called the Wasserstein distance.
\end{dfn}

The following theorem is known regarding the Wasserstein distance :
\begin{thm}
	For any $p \in [1, \infty)$,
	$W_p$ is a metric on $P_p(X)$,
	meaning that it satisfies the following three axioms of a metric,
	where $\mu, \nu, \lambda \in P_p(X)$ :
	\begin{enumerate}
		\item $W_p(\mu, \nu) \geq 0$ and $W_p(\mu, \nu) = 0$ if and only if $\mu = \nu$.
		\item $W_p(\mu, \nu) = W_p(\nu, \mu)$.
		\item $W_p(\mu, \lambda) \leq W_p(\mu, \nu) + W_p(\nu, \lambda)$.
	\end{enumerate}
\end{thm}

In this paper,
we refer to the Wasserstein distance as the optimal transport distance.
The optimal transport between discrete probability distributions can be understood
as the optimal transport problem for Dirac measures $\delta_x\ (x \in X)$.
Using vectors and matrices, it can be expressed as follows :
\begin{dfn}
	Let $\bm{p} = (p_i)_{i=1}^{N}, \bm{q} = (q_j)_{j=1}^{N}$.
	The optimization problem
	\begin{gather*}
		\textrm{OT}(\bm{p}, \bm{q}, C) := \min_{P \in \mathbb{R}^{N \times N}}\left\{\sum_{i=1}^{N}\sum_{j=1}^{N}C_{ij}P_{ij}\ \middle|\ \sum_{j=1}^{N}P_{ij} = p_i, \sum_{i=1}^{N}P_{ij} = q_j\right\}
	\end{gather*}
	to determin the transport matrix $P=(P_{ij})_{i,j=1}^{N}\ (P_{ij} \geq 0)$
	for the cost matrix $C=(C_{ij})_{i,j=1}^{N}$ is called the optimal transport problem.
	In addition, the solution of this optimization problem is called the optimal transport distance.
\end{dfn}

What we aimed to investigate was the changes in clustering patterns in the partially ordered phase.
To evaluate the variation in clustering patterns,
we can transform the clustering pattern at each step into a distribution that contains information about the number of elements in each cluster.
This allows us to use the optimal transport distance for evaluation.
For this purpose,
we consider the cluster distribution,
which is defined through a mapping whose domain is the set of all clustering patterns :
\begin{dfn}
	For any $N \in \mathbb{N}, k \leq N$,
	let
	\begin{gather*}
		\Omega_k := \{(N_1,N_2,\ldots,N_k) \in \mathbb{N}^k \mid N_1+N_2+\cdots+N_k=N\}
	\end{gather*}
	be the set of all clustering patterns of degree $k$,
	and
	\begin{gather*}
		\Sigma_{N} := \left\{(p_1,p_2,\ldots,p_N)\ \middle|\ \sum_{i=1}^Np_i=1\right\}.
	\end{gather*}
	We define a map $T : \Omega_k \to \Sigma_N$ as follows:
	\begin{gather*}
		T(N_1,N_2,\ldots,N_k) = (p_1,p_2,\ldots,p_N),
		\quad
		p_i = \frac{1}{N}\sum_{\substack{j=1,2,\ldots,N;\\ \textrm{}N_j=i}}N_j.
		\quad
		i=1,2,\ldots,N
	\end{gather*}
\end{dfn}
\begin{rem}
	This definition is formulated such that the optimal transport distance between the state in the coherent phase and the state in the turbulent phase is maximized compared to other cases
	(see Example \ref{exam:clus_dist}).
\end{rem}

\begin{exam}
	\label{exam:clus_dist}
	Fig.~\ref{fig:exam_clus_dist} shows examples of the cluster distributions for some clustering patterns.
	For simplicity, we assume $N=10$ here,
	and we abbreviate an clustering pattern,
	using $431^3$ in place of $(4, 3, 1, 1, 1)$,
	for example.
	\begin{enumerate}[label=\textbf{(\alph*)}]
		\item In the case of the clustering pattern is $(10)$,
		i.e., all elements are synchronized, the cluster distribution is
		$(\underbrace{0,\ldots,0}_{9},1)$:
		\begin{gather*}
			p_{10} = \frac{N_1}{N} = \frac{10}{10} = 1,
			\quad
			p_i = 0.
			\quad
			(i = 1,2,\ldots,9)
		\end{gather*}
		\item In the case of the clustering pattern is $1^{10}$,
		i.e., all elements are out of synchronization, the cluster distribution is
		$(1,\underbrace{0,\ldots,0}_{9})$:
		\begin{gather*}
			p_1
			= \sum_{j=1,2,\ldots,10;N_j=1}\frac{N_j}{10}
			= 10 \times \frac{1}{10}
			= 1,
			\quad
			p_i = 0.
			\quad
			(i = 2,3,\ldots, 10)
		\end{gather*}
		\item In the case of the clustering pattern is $5 4 1$,
        i.e., 3-clusters,
		the cluster distribution is
		$(0.1,0,0,0.4,0.5,\underbrace{0,\ldots,0}_{5})$:
		\begin{gather*}
			p_1 = \sum_{j=1,2,3;N_j=1}\frac{N_j}{10}
			= \frac{1}{10},
			\quad
			p_4 = \sum_{j=1,2,3;N_j=4}\frac{N_j}{10}
			= \frac{4}{10},\\
			p_5 = \sum_{j=1,2,3;N_j=5}\frac{N_j}{10}
			= \frac{5}{10},
			\quad
			p_i = 0.
			\quad
			(i = 2,3,6,\ldots,10)
		\end{gather*}
	\end{enumerate}

	\begin{figure}[t]
		\centering
		\begin{minipage}[t]{0.3\columnwidth}
			\centering
			\includegraphics[width=\columnwidth]{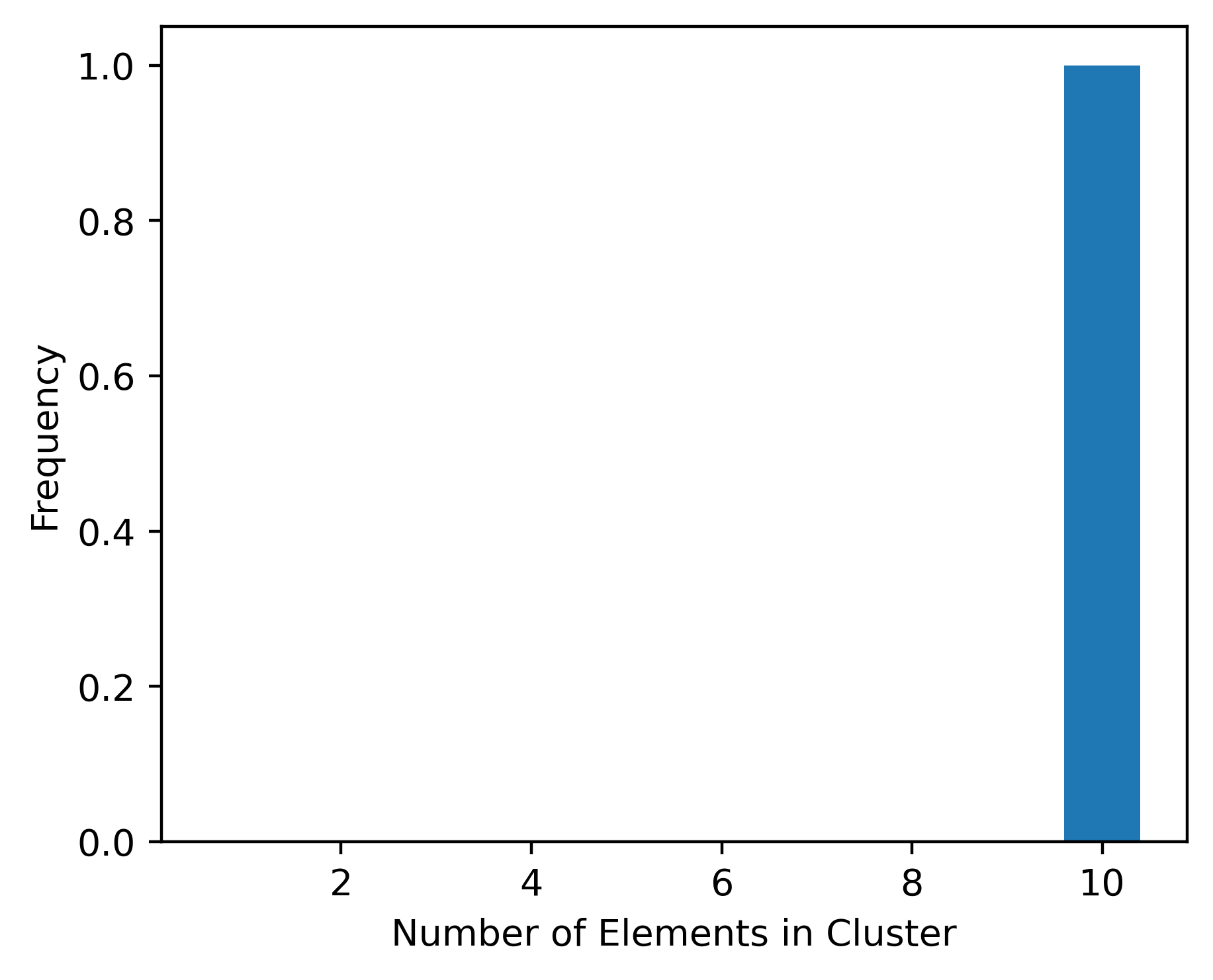}
			\subcaption{$(10) = 10$}
      \label{fig:4-a}
		\end{minipage}
		\hspace{0.02\columnwidth}
		\begin{minipage}[t]{0.3\columnwidth}
			\centering
			\includegraphics[width=\columnwidth]{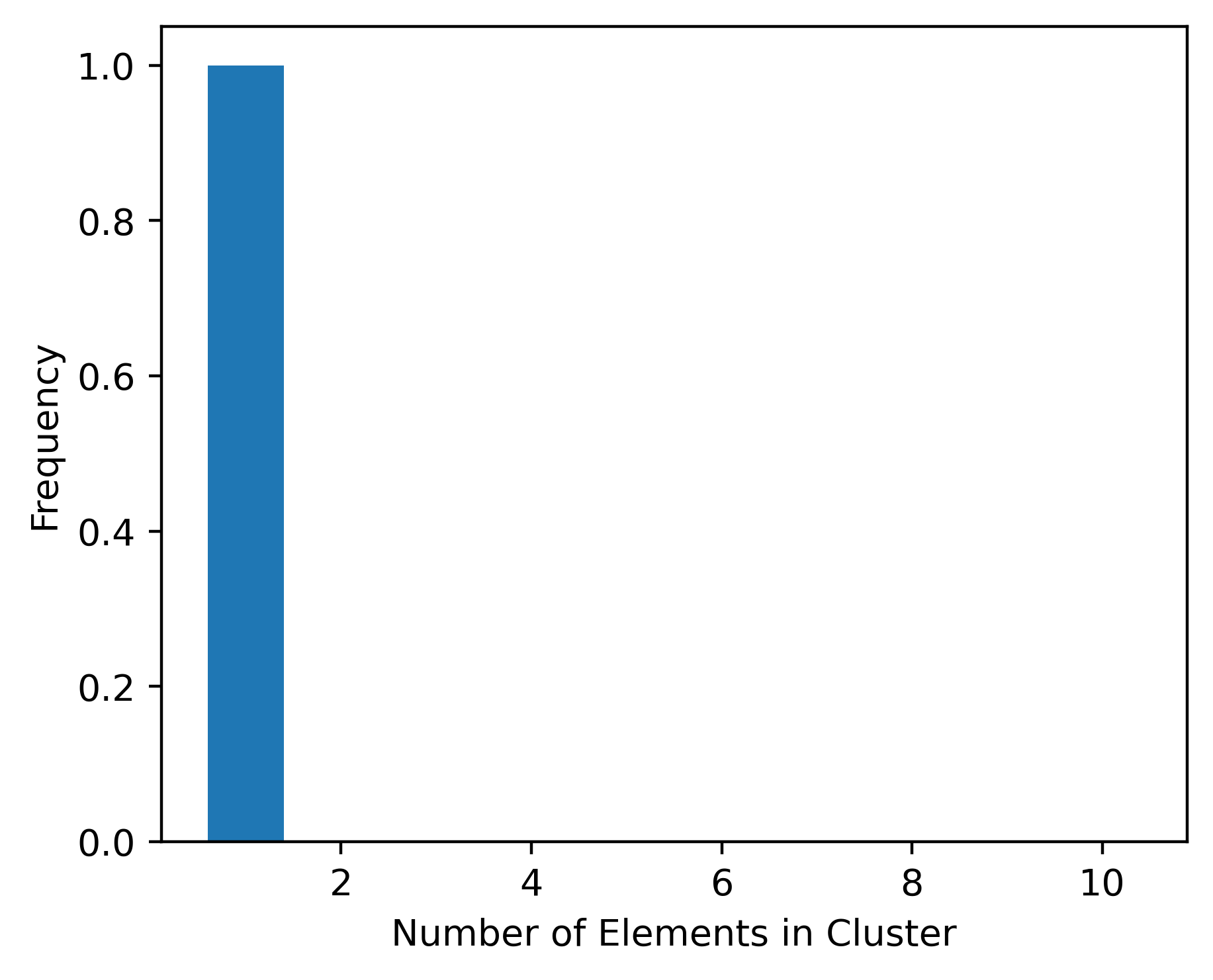}
			\subcaption{$(\underbrace{1,\ldots,1}_{10}) = 1^{10}$}
      \label{fig:4-b}
		\end{minipage}
		\hspace{0.02\columnwidth}
		\begin{minipage}[t]{0.3\columnwidth}
			\centering
			\includegraphics[width=\columnwidth]{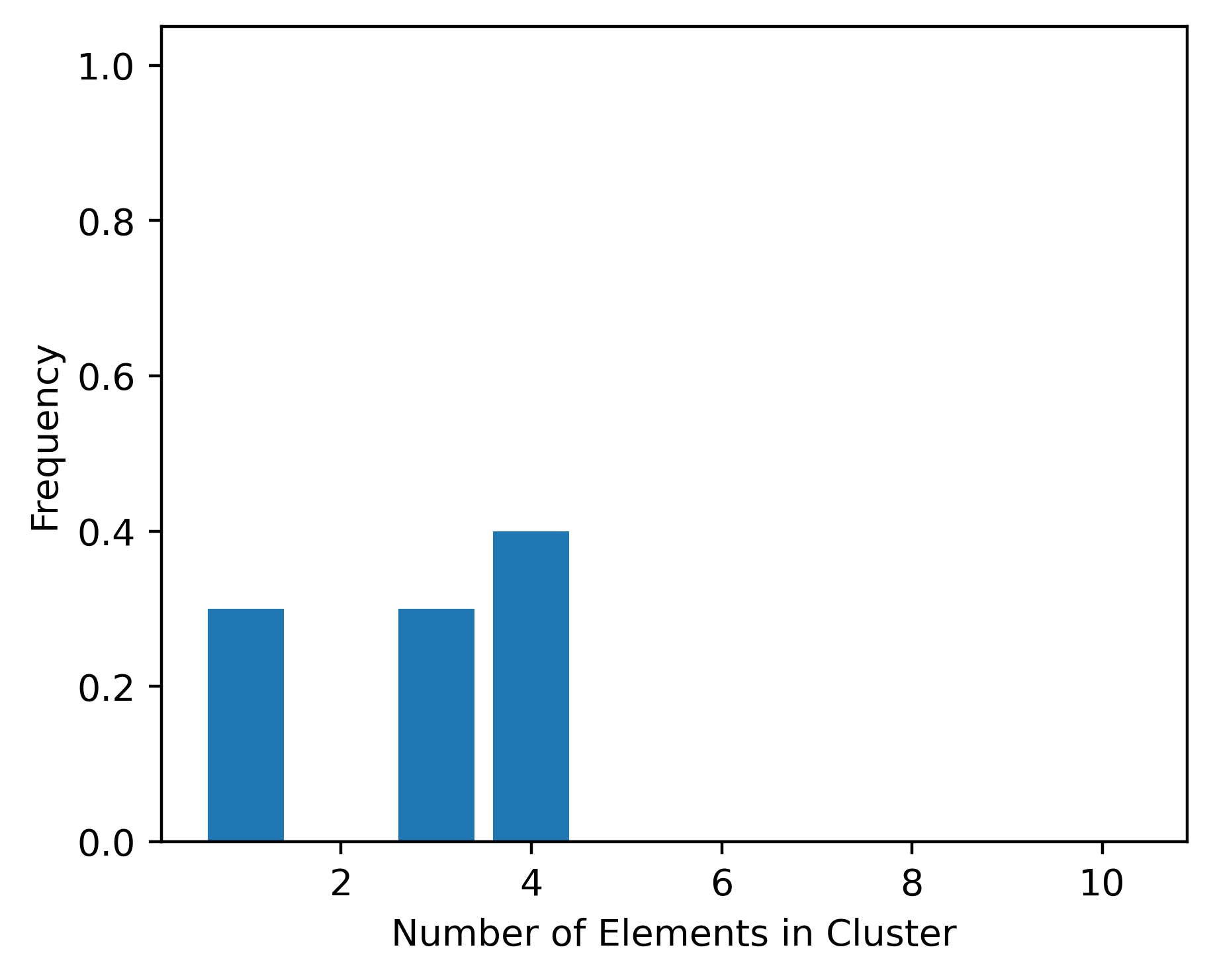}
			\subcaption{$(4, 3, 1, 1, 1) = 4 3 1^3$}
		\end{minipage}
		\caption{The examples of the cluster distributions for each the clustering patterns $(N_1,N_2, \ldots, N_k)$.
		Note that it has been normalized so that the total sum equials one.
		Fig.~\ref{fig:4-a} and Fig.~\ref{fig:4-b} are the states in the coherent phase and the turbulent phase, respectively.}
		\label{fig:exam_clus_dist}
	\end{figure}
\end{exam}

\section{Main results}
\subsection{Time series of optimal transport distance and their time averages}
\begin{figure}[t]
	\centering
	\begin{minipage}[t]{0.45\columnwidth}
		\centering
		\includegraphics[width=0.9\columnwidth]{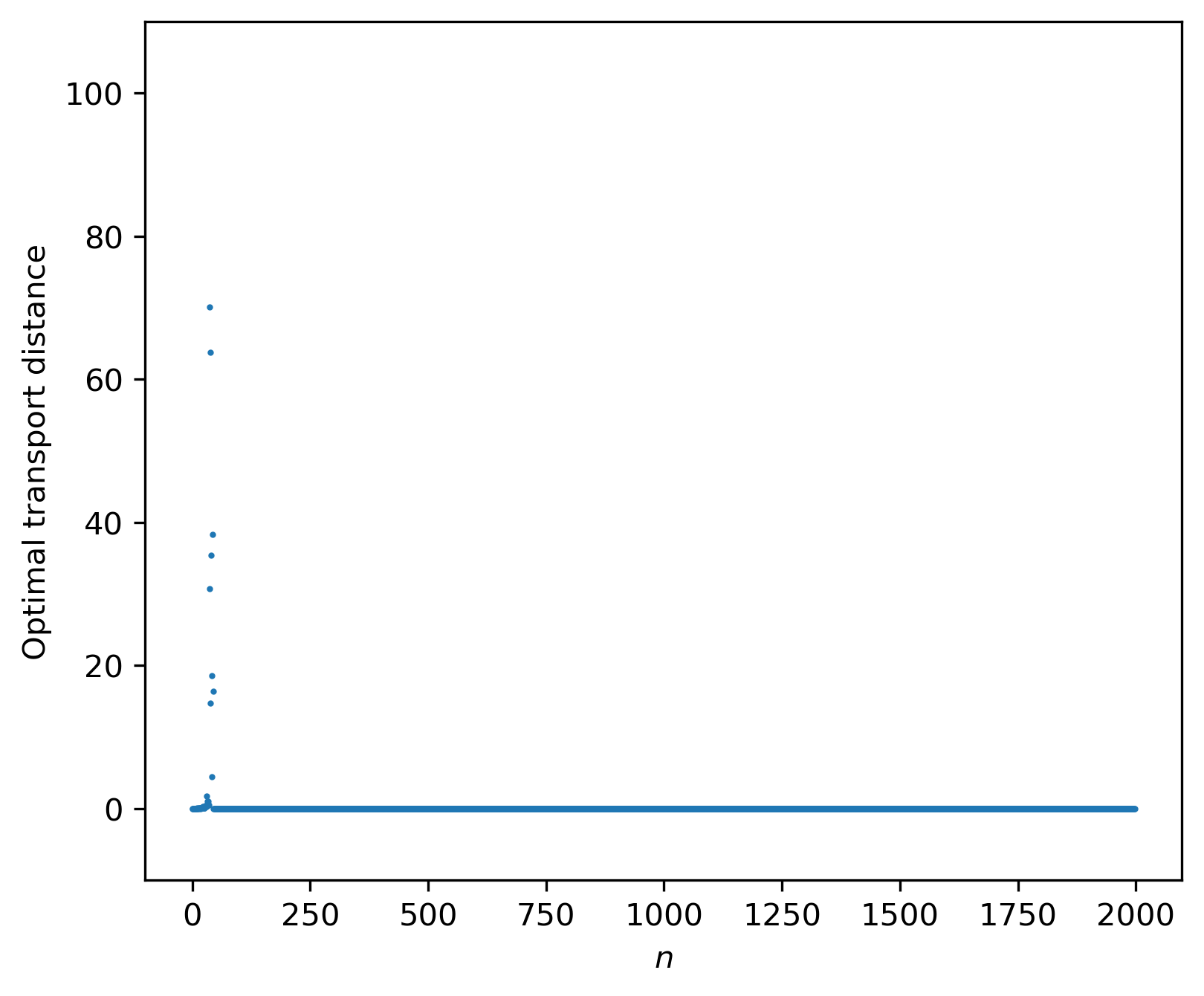}
		\subcaption{$\varepsilon = 0.5$ (Coherent phase)}
		\label{fig:3-a}
	\end{minipage}
	\hspace{0.05\columnwidth}
	\begin{minipage}[t]{0.45\columnwidth}
		\centering
		\includegraphics[width=0.9\columnwidth]{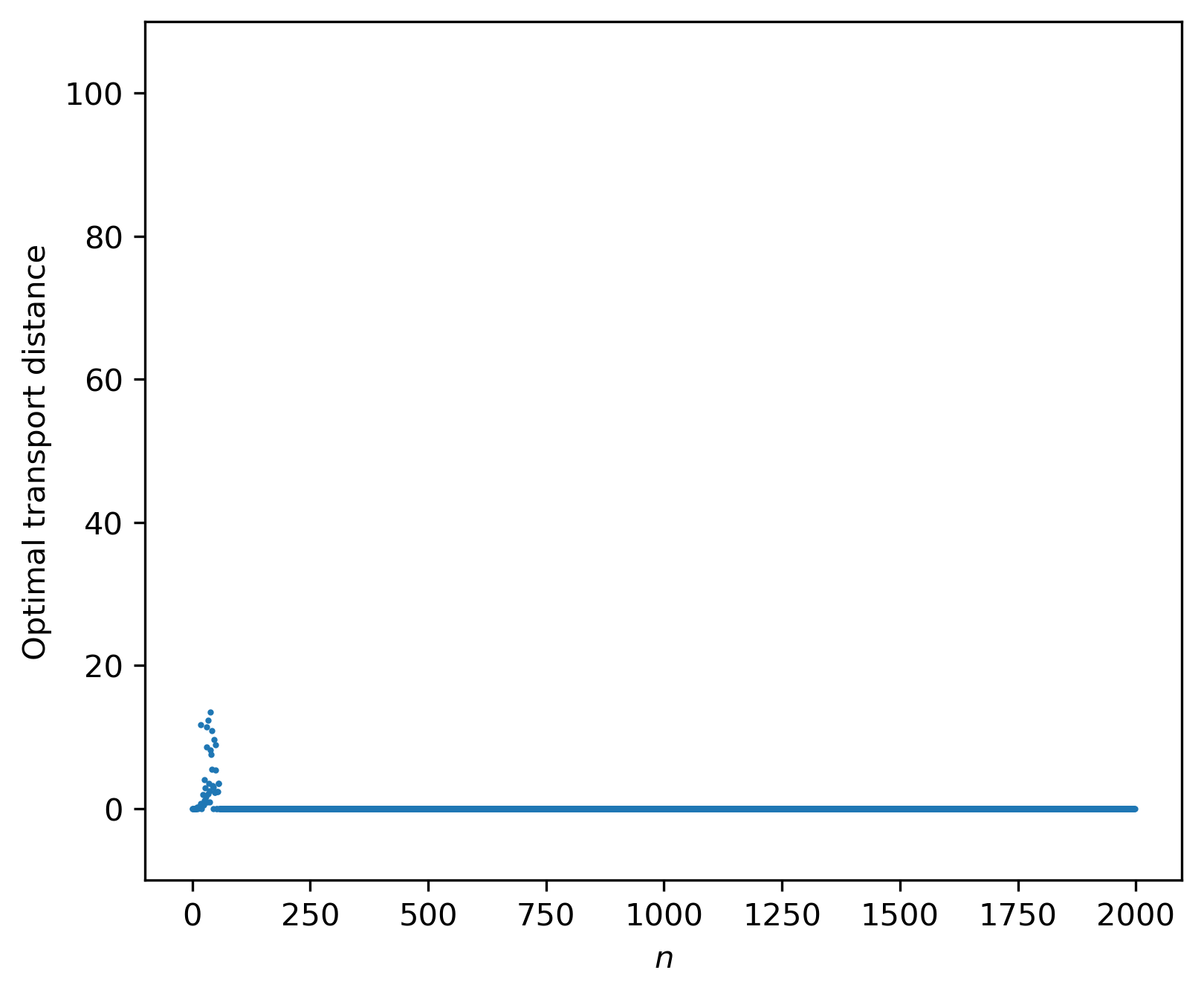}
		\subcaption{$\varepsilon = 0.2$ (Ordered phase)}
		\label{fig:3-b}
	\end{minipage}
	\vspace{10pt}\\
	\begin{minipage}[b]{0.45\columnwidth}
		\centering
		\includegraphics[width=0.9\columnwidth]{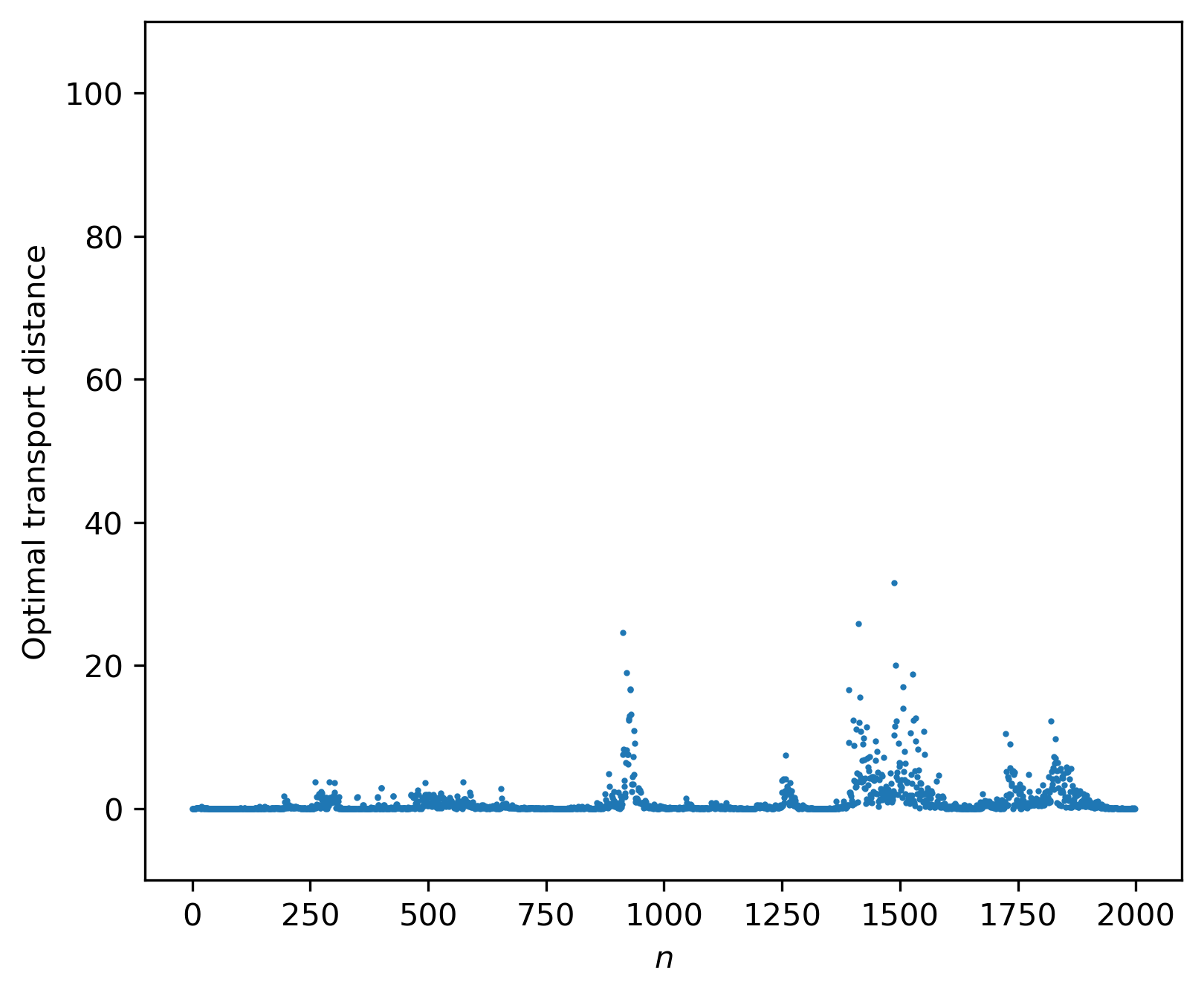}
		\subcaption{$\varepsilon = 0.3$ (Partially ordered phase)}
		\label{fig:3-c}
	\end{minipage}
	\hspace{0.05\columnwidth}
	\begin{minipage}[b]{0.45\columnwidth}
		\centering
		\includegraphics[width=0.9\columnwidth]{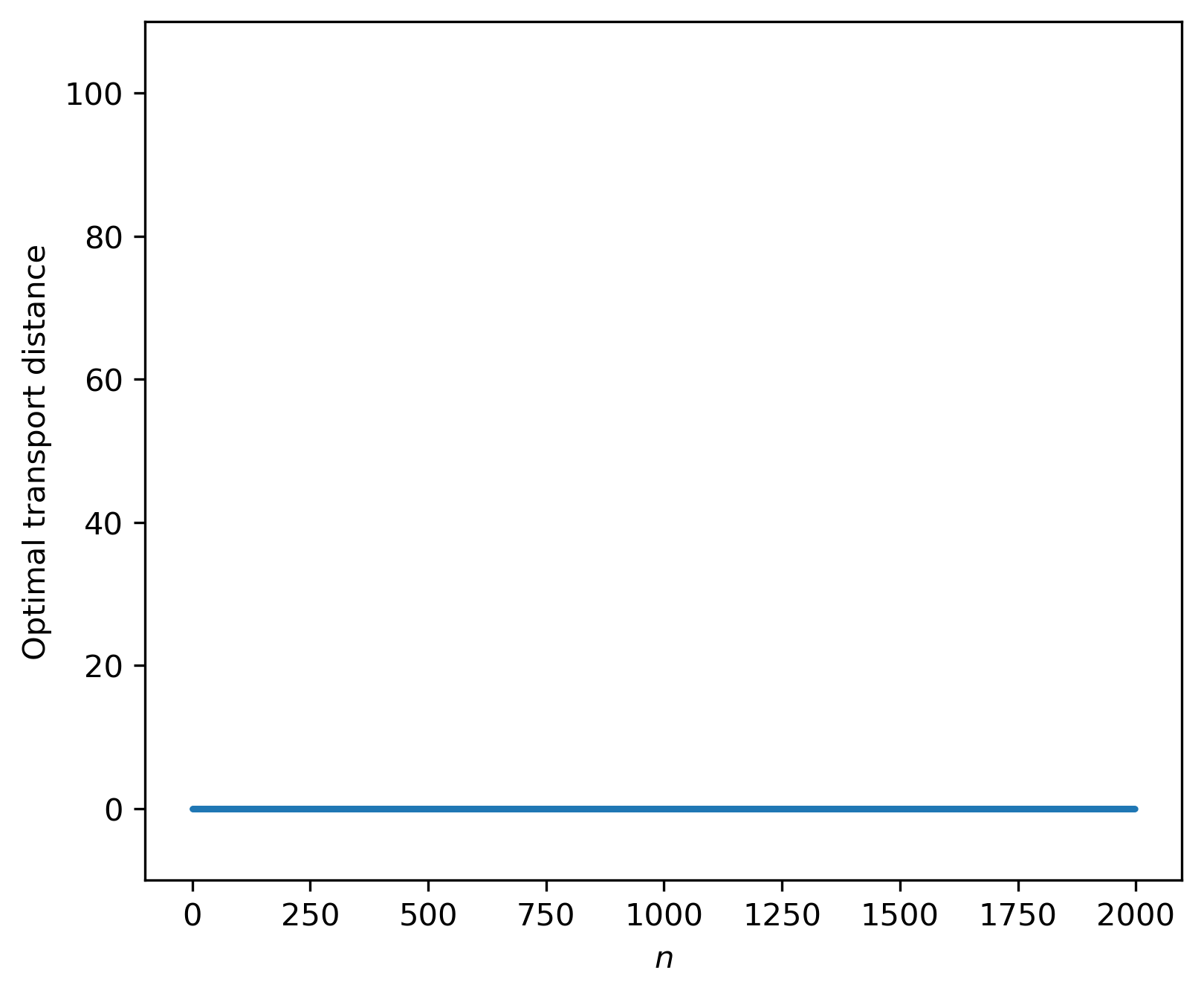}
		\subcaption{$\varepsilon = 0.0$ (Turbulent phase)}
		\label{fig:3-d}
	\end{minipage}
	\caption{Examples of the optimal transport distance for each $\varepsilon$.
	Let $\alpha = 3.8, N = 100, \delta = 10^{-6}$ and the initial conditions are identical for any $\varepsilon$.}
	\label{fig:3}
\end{figure}

\begin{figure}[t]
	\centering
  \begin{minipage}[b]{0.48\columnwidth}
		\centering
		\includegraphics[width=\columnwidth]{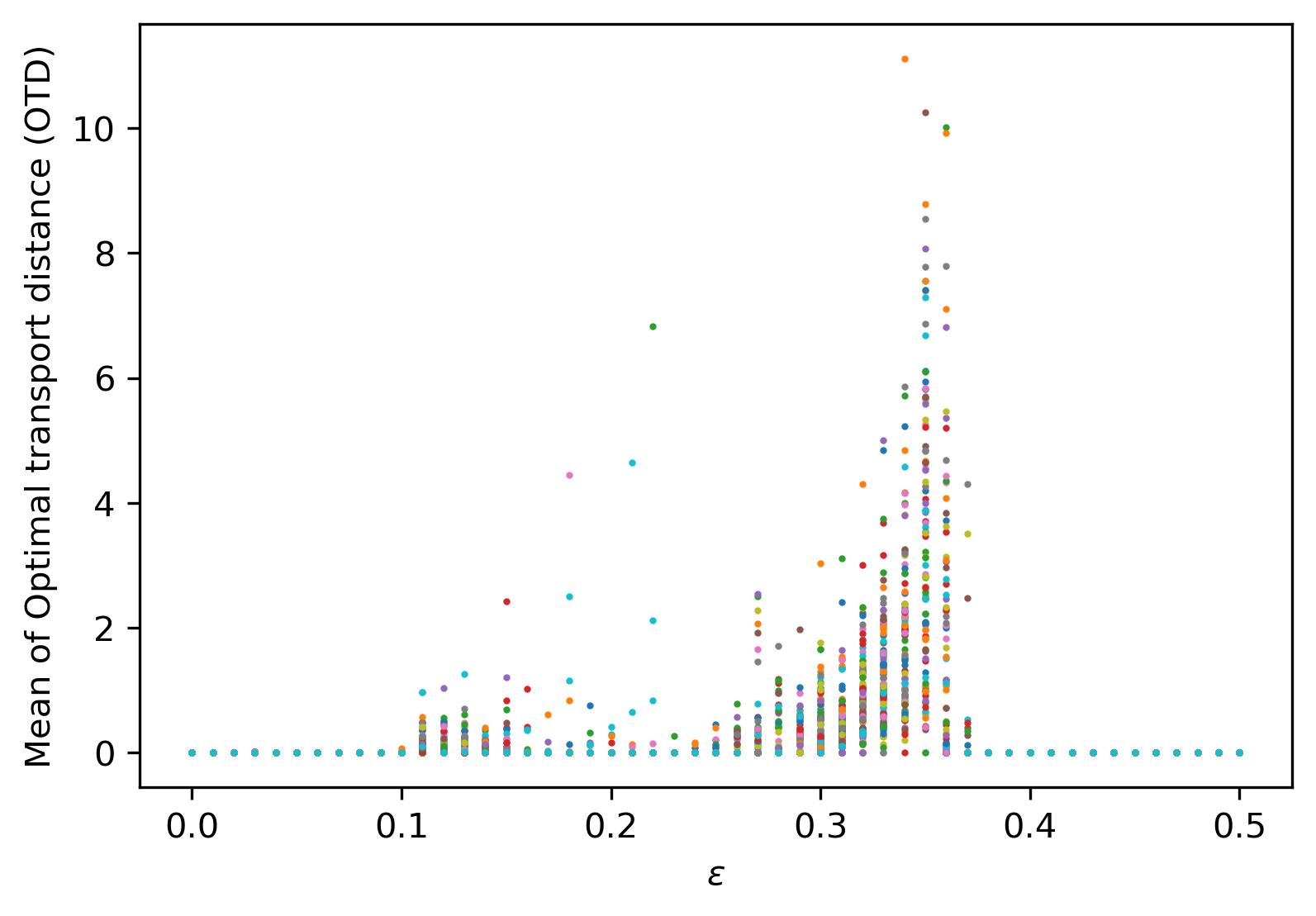}
		\caption{The time average of optimal transport distance.
		Here, let $N=100$, $\alpha = 3.8$.}
		\label{fig:4}
	\end{minipage}
	\hfill
  \begin{minipage}[b]{0.48\columnwidth}
		\centering
		\includegraphics[width=\columnwidth]{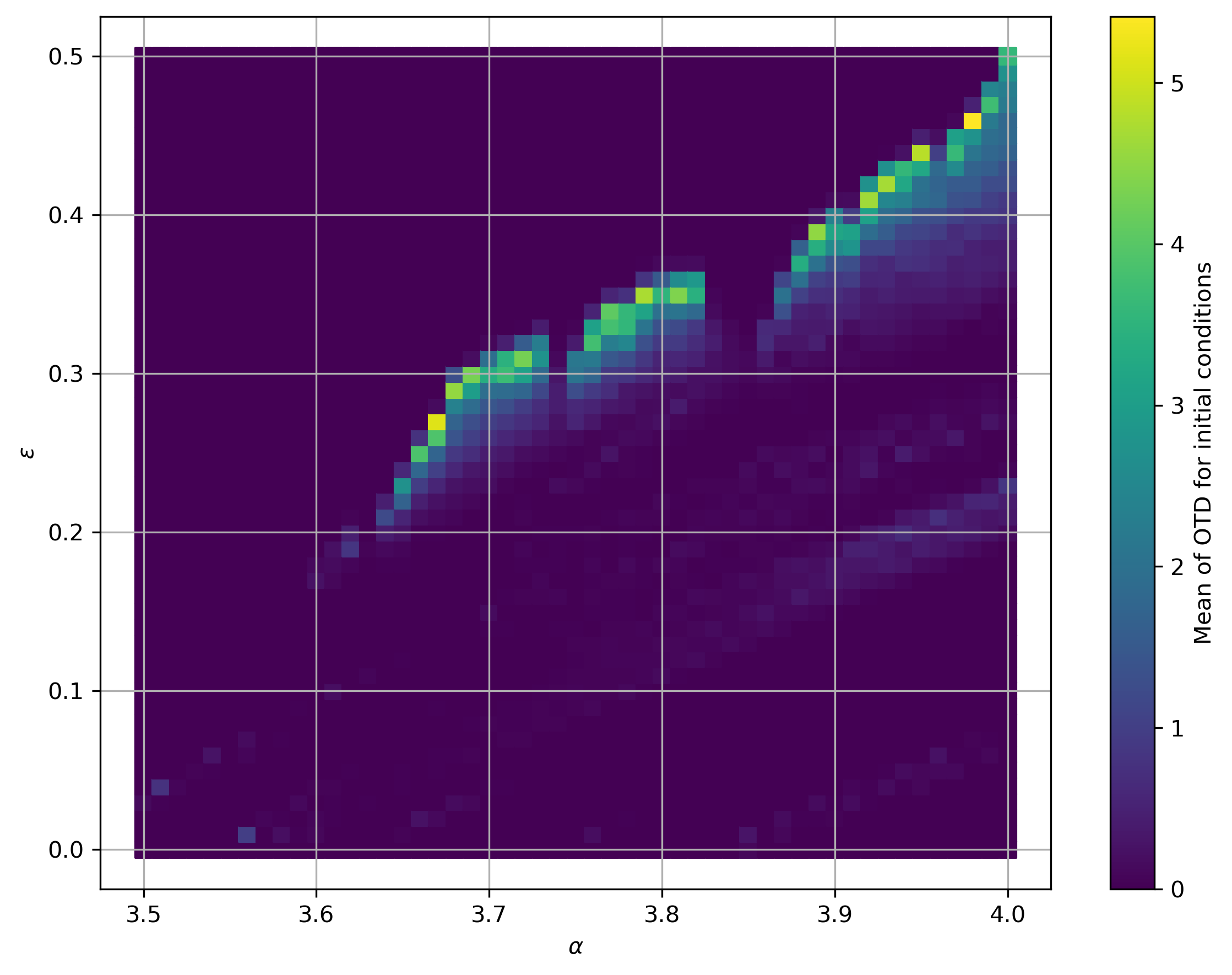}
		\caption{The average optimal transport distance for 100 initial conditions.
		Here, let $N = 100$.}
		\label{fig:5}
	\end{minipage}
\end{figure}

We observed intermittent phenomena in the partially ordered phase (Fig.~\ref{fig:3-c})
and obtained a higher time averages of the time series of optimal transport distance (Fig.~\ref{fig:5}).

Fig.~\ref{fig:3} shows the time series of the optimal transport distance for each $\varepsilon$ when $N=100, \alpha = 3.8$ is fixed.
In the parameter region corresponding to the coherent phase, the ordered phase, and the turbulent phase,
the clustering pattern converges after a transient process and remains unchanged.
As a result, the time series of the optimal transport distance also becomes zero after a certain time,
which can be confirmed from the analysis (Fig.~\ref{fig:3-a}, \ref{fig:3-b}, \ref{fig:3-d}).
In contrast, in the parameter region corresponding to the partially ordered phase,
the clustering pattern does not converge.
This characteristic is reflected in the time series of the optimal transport distance,
which continues to fluctuate over time (Fig.~\ref{fig:3-c}).

Fig.~\ref{fig:4} shows the time-averaged optimal transport distance for each $\varepsilon$ when $\alpha = 3.8$ is fixed.
The results for 100 different initial conditions are overlaid,
and the average is computed using the last 1,000 steps of a time series of length 2,000,
with the first 1,000 steps treated as the transient period.
In the parameter region corresponding to the coherent phase and  the turbulent phase,
the clustering patterns converge regardless of the initial conditions.
As a result, the time-averaged optimal transport distance,
which quantifies the variation in clustering patterns,
becomes zero.
In the parameter region corresponding to the ordered phase,
the clustering pattern may depend on the initial conditions,
leading to different final states.
However, since the pattern remains unchanged after a certain time,
the time-averaged optimal transport distance also becomes zero in this region.
In contrast,
in the partially ordered phase,
the time-averaged optimal transport distance remains greater than zero,
indicating persistent fluctuations in the clustering pattern.
In particular, it is observed that the time-averaged optimal transport distance tends to be larger in the partially ordered phase I\hspace{-1.2pt}I compared to the partially ordered phase I.
Fig.~\ref{fig:5} presents the results of averaging the time-averaged optimal transport distance over 1,000 initial conditions for each parameter.
It can be observed that the values tend to be higher in the parameter region corresponding to the partially ordered phase I\hspace{-1.2pt}I in the phase diagram of GCM (see Fig.~\ref{fig:phase_diagram}).

\subsection{Strength of attractor-ruins}
We evaluated the strength of attractor-ruins using an optimal transport distance and the effective dimension.
As a result, the strength of attractor-ruins is higher in the partially ordered phase.

\begin{figure}[t]
	\centering
	\begin{minipage}[t]{0.45\columnwidth}
		\centering
		\includegraphics[width=0.9\columnwidth]{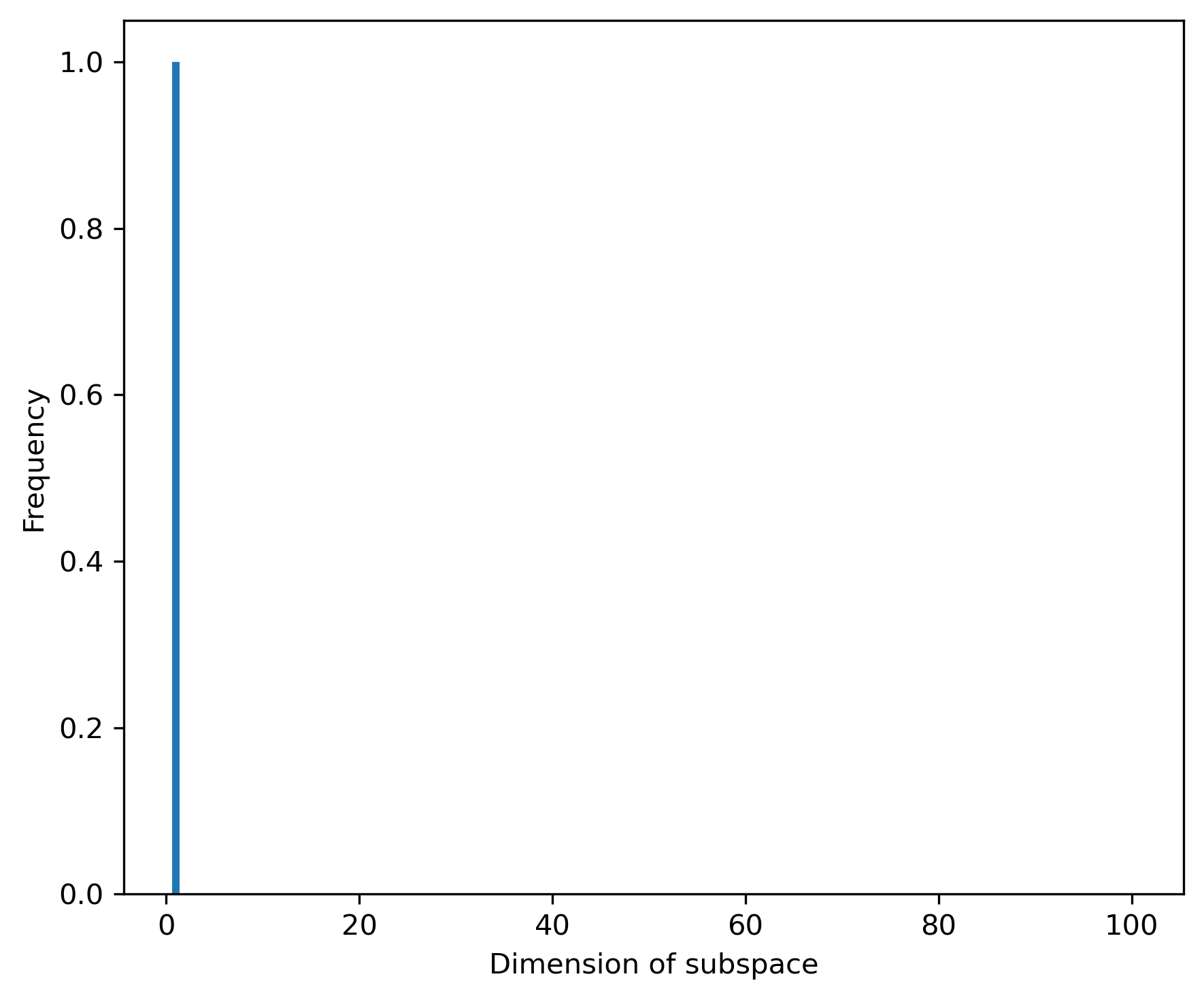}
		\subcaption{$(\alpha, \varepsilon) = (3.8, 0.5)$.}
		\label{fig:example_entropy_dist_1}
	\end{minipage}
	\hspace{0.05\columnwidth}
	\begin{minipage}[t]{0.45\columnwidth}
		\centering
		\includegraphics[width=0.9\columnwidth]{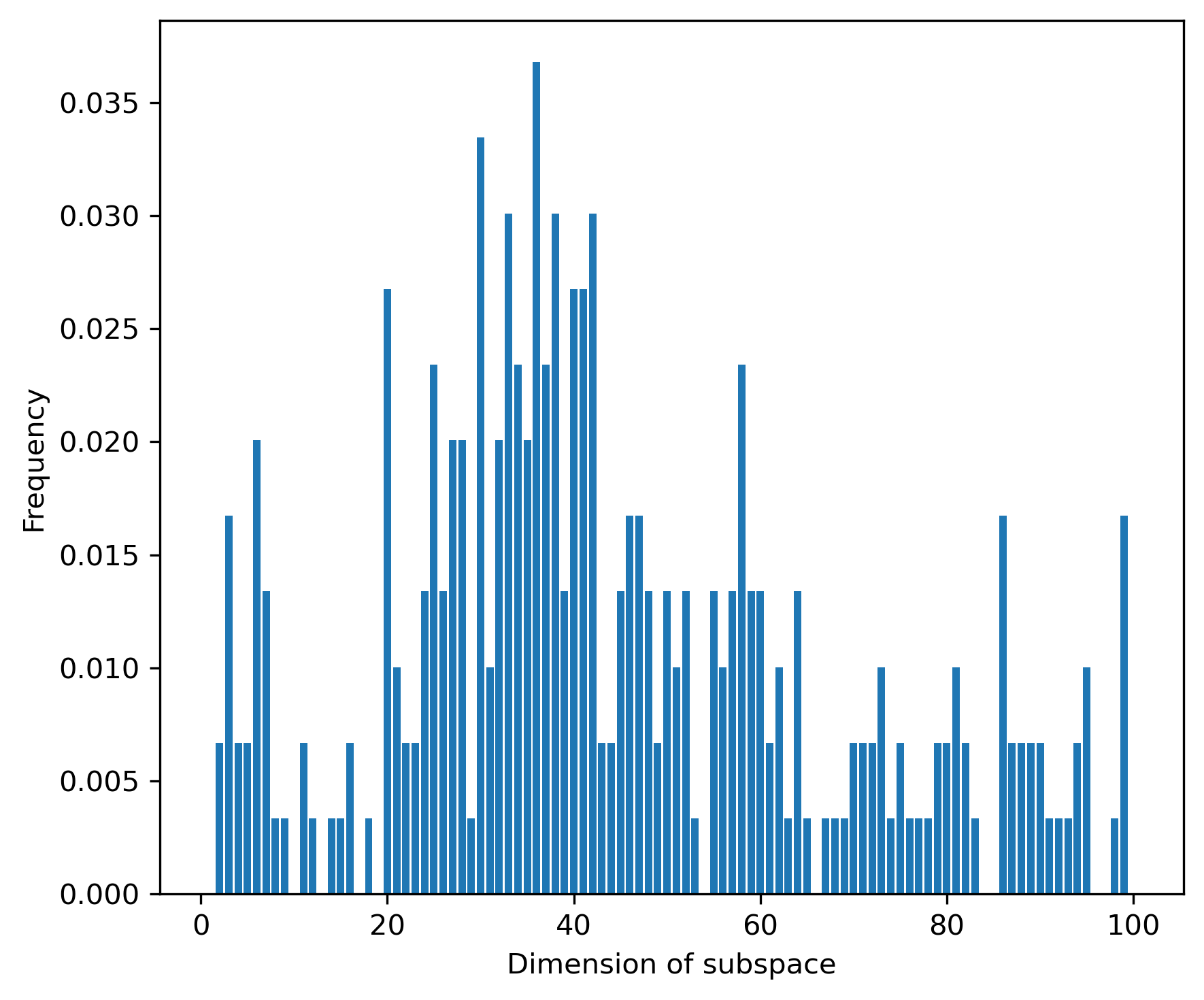}
		\subcaption{$(\alpha, \varepsilon) = (3.8, 0.33)$.}
		\label{fig:example_entropy_dist_2}
	\end{minipage}
	\caption{Example of discrete distribution for some initial condition and different parameters.
	Fig.~8(a) corresponds to the coherent phase,
	while Fig.~8(b) corresponds to the partially ordered phase I\hspace{-1.2pt}I.}
	\label{fig:example_entropy_dist}
\end{figure}

\begin{figure}[t]
	\centering
  \begin{minipage}[b]{0.48\columnwidth}
		\centering
		\includegraphics[width=\columnwidth]{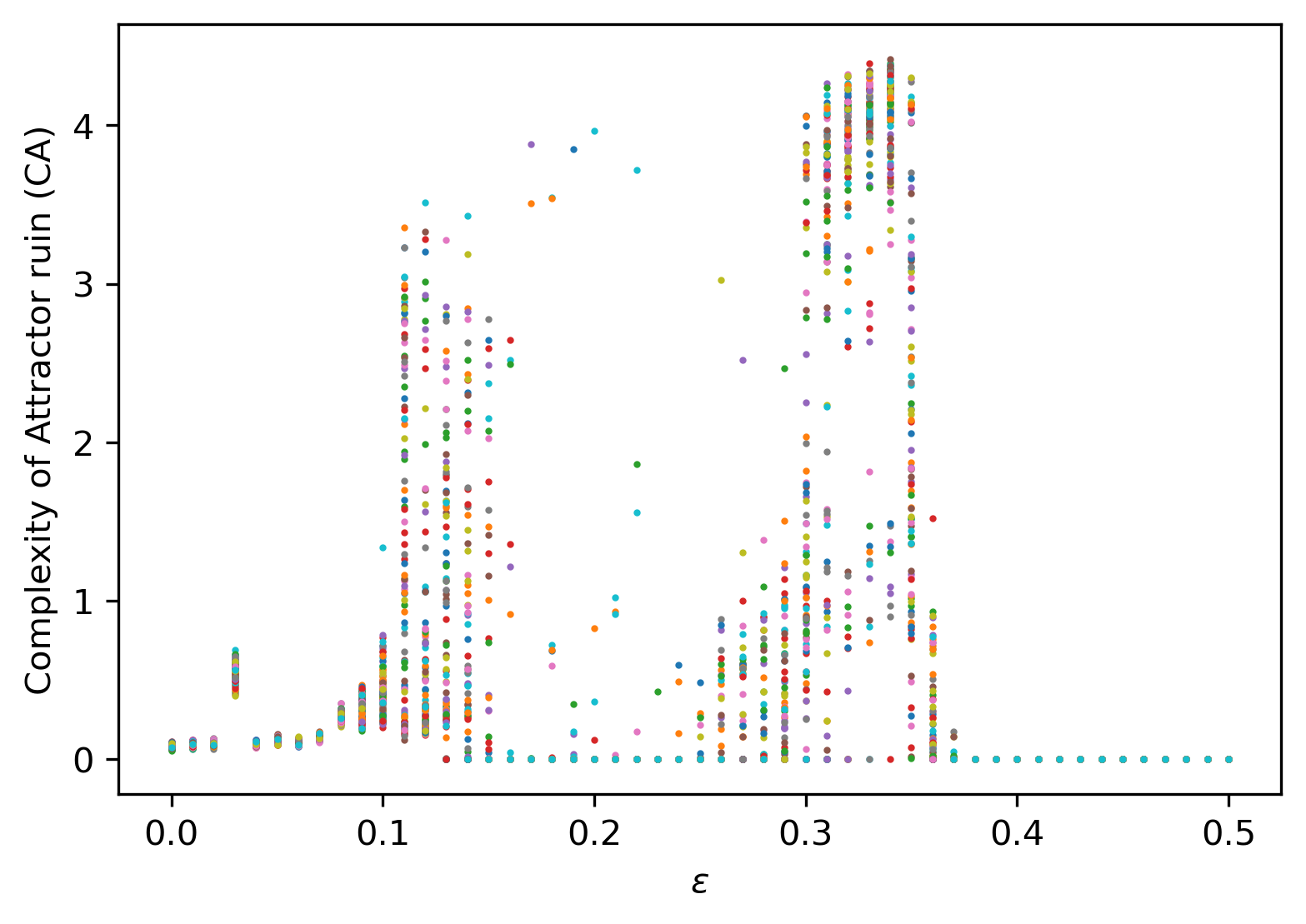}
    \caption{The strength of the attractor-ruins.
		Here, let $N=100$, $\alpha = 3.8$.}
		\label{fig:6}
	\end{minipage}
	\hfill
	\begin{minipage}[b]{0.48\columnwidth}
		\centering
		\includegraphics[width=\columnwidth]{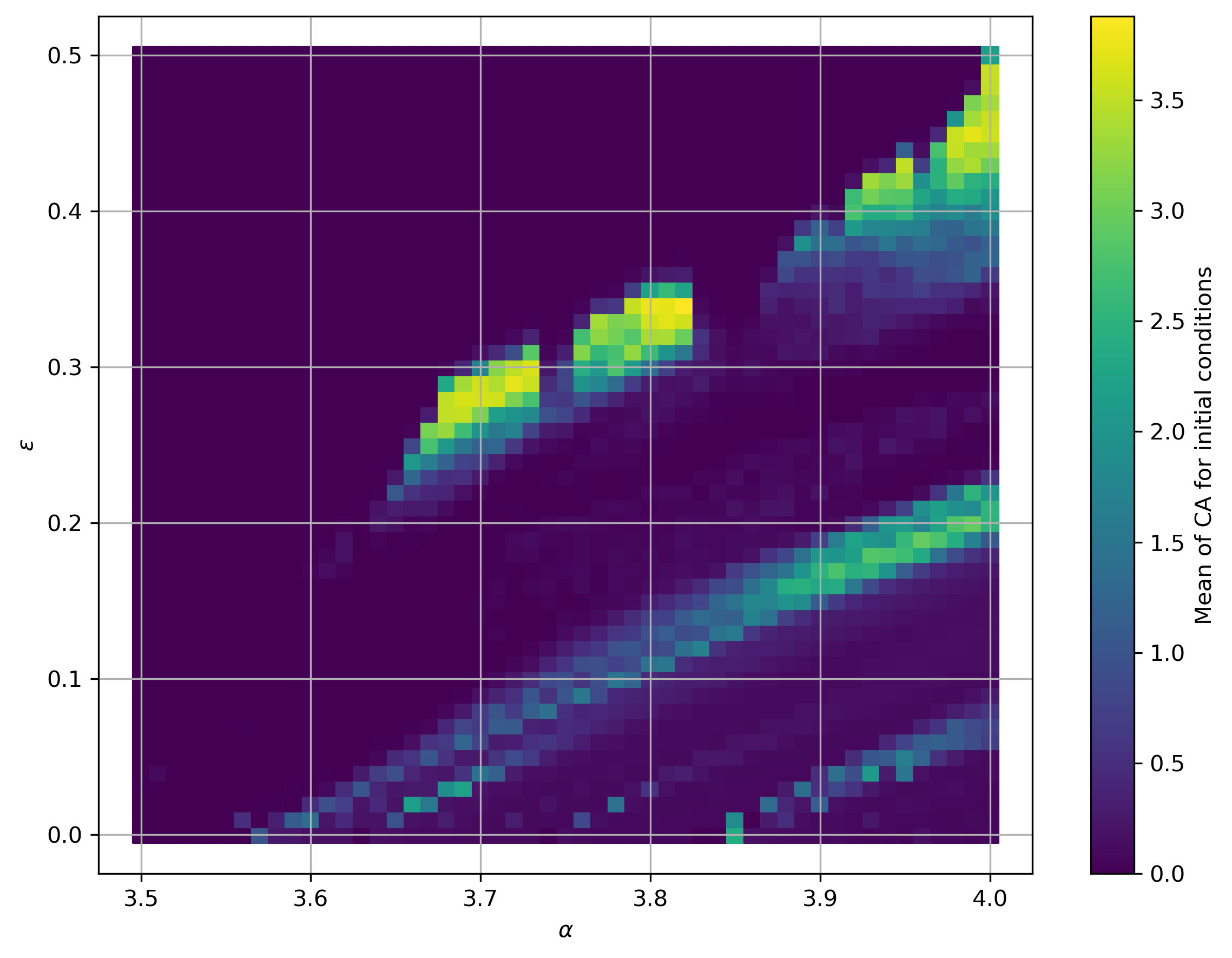}
    \caption{The strength of the attractor-ruins.
		Here, let $N=100$.
		This results are averages for 100 initial conditions.}
		\label{fig:7}
	\end{minipage}
\end{figure}

Fig.~\ref{fig:6} shows the strength of the attractor-ruins calculated using Shannon entropy.
Here, we assume $\alpha = 3.8$, $N=100$, $\delta = 10^{-6}$,
and the results for 100 initial conditions are overlaid,
and the average is computed using the last 10,000 steps of a time series of length 11,000,
with the first 1,000 steps treated as the transient period.
Furthermore, the strength of the attractor-ruins in this context refers to the results obtained through the following calculation:
At first, the time series of effective dimension and optimal transport distance are computed from the time series of orbits.
The effective dimension corresponding to the time points where the optimal transport distance is zero is then extracted,
and the distribution of these data points is obtained.
Shannon entropy is calculated as a feature of this discrete distribution (see Fig.~\ref{fig:example_entropy_dist}).
We consider the Shannon entropy computed in this way to be the strength of the attractor-ruins.
First, in the coherent phase,
the system converges to a single-cluster state regardless of the initial conditions,
resulting in a Shannon entropy of zero for any initial condition (see Fig.~\ref{fig:example_entropy_dist_1} as an example).
Similarly, in the turbulent phase,
the effective dimension approaches the number of elements irrespective of the initial conditions,
leading to a low Shannon entropy.
In the ordered phase,
for most initial conditions,
the clustering pattern converges in the same manner as in the coherent phase,
yielding a Shannon entropy of zero.
Finally, in the partially ordered phase,
the effective dimension at the steps where the optimal transport distance becomes zero exhibits a wide range of values (see Fig.~\ref{fig:example_entropy_dist_2}),
resulting in a high Shannon entropy that is highly dependent on the initial conditions.

Fig.~\ref{fig:7} presents the results of averaging the Shannon entropy over 100 initial conditions for each parameter.
It can observed that the values tend to be higher in the partially ordered phase (see Fig.~\ref{fig:phase_diagram}).

\section{Summary and Discussion}
In this study, we have characterized each region of parameters for GCM from the perspective of clustering.
First, we evaluated the amount of change in the clustering pattern using the optimal transport distance.
As a result, we show that the time series of the optimal transport distance is zero for phases in which the clustering pattern does not change,
such as the coherent phase,
while the time series is higher for the partially ordered phase in which the clustering pattern changes (see Fig.~\ref{fig:3}-\ref{fig:5}).
In particular,
since the optimal transport distance is observed to be higher in the partially ordered phase I\hspace{-1.2pt}I,
it can be considered that this measure captures the characteristic features emerging in this phase.
From the perspective of clustering,
the partially ordered phase I\hspace{-1.2pt}I can be interpreted as significant variations in clustering patterns.
Second, we evaluated the strength of the attractor-ruins using the Shannon entropy derived from the distribution,
which can be obtained by combining the time series of the effective dimension and the optimal transport distance (see Fig.~\ref{fig:6}-\ref{fig:7}).
The results show that the Shannon entropy is higher,
i.e., the strength of attractor-ruins is higher,
in the partially ordered phase,
indicating that attractor-ruins with various dimensions attract orbits.

By applying the methodology presented in this paper,
we can provide a novel classification of phase in GCM,
focusing on the fluctuations of clustering patterns.
In future works, by using optimal transport distance,
we would like to investigate the relationship with the Basin structure as a Milner attractor
and the statistics when viewed as an intermittent phenomenon.

\appendix
\section*{Appendix}\label{sec:appendix}
A problem known as ``apparent attraction'' occurs when a pair of elements, $x(i)$ and $x(j)$,
are not truly equal but are computationally judged to be equal due to the limitations of the number of digits the computer can handle.
To resolve this issue,
we introduce a small amount of noise at each step when calculating the orbit.
If the system is chaotic,
this small error induced by the noise is expected to grow,
causing the values of the two elements to diverge.
In the results of this paper,
the same noise is applied and calculations are performed for all initial values and all time points.

The definition of effective dimension is provided in Def. \ref{dfn:effdim}.
According to this definition, the effective dimension involves searching for all subsets $H_\sigma$ where the intersection with the $\delta$-neighborhood of point $x$ is not empty,
and then identifying the subset with the smallest dimension.
However, the number of such subsets $H_\sigma$ increases exponentially with the system dimension $N$,
which could lead to an enormous computational cost.
Therefore, in this paper,
we consider the ``effective precision cluster,''
which can be obtained by determining the synchronization or asynchrony between elements with respect to
\begin{gather*}
	\bar{x}_\delta(i) = \delta \times \left[\frac{x(i)}{\delta}\right].
\end{gather*}
It should be noted that the effective dimension in the following refers to the number of effective precision in this sense.

\section*{CRediT authorship contribution statement}
Koji Wada: Conceptualization, Methodology, Software, Formal analysis, Investigation, Writing - Original Draft preparation, Funding acquisition.
Takao Namiki: Supervision, Writing-review and editing, Funding acquisition.

\section*{Acknowledments}
This work was supported by JST SPRING, Grant Number JPMJSP2119 (KW)
and JSPS KAKENHI Grant Number 23K2578503 (TN).

\bibliographystyle{jabbrv} 
\bibliography{arx_refs}

\end{document}